\documentclass[11pt]{article}
\usepackage{amsmath}
\usepackage[dvips]{graphicx}
\usepackage{amssymb, amsmath, amsfonts, latexsym}
\usepackage{amssymb}
\usepackage{color}

\newcommand{\qed}{\hfill \rule{1.6mm}{1.6mm}}
\newcommand{\no}{\nonumber}
\newcommand{\noi}{\noindent}

\catcode`\@=11

\def\theequation{\@arabic{\c@section}.\@arabic{\c@equation}}
\catcode`\@=12
\numberwithin{equation}{section}
\title{ Sufficient stochastic maximum principle for the optimal control of semi-Markov modulated jump-diffusion with application to Financial optimization.}
\author{Amogh Deshpande  \footnote {Department of Statistics, University of Warwick, Coventry,  CV47AL,UK.
Email: addeshpa@gmail.com}}
\date{}
\begin{document}
\maketitle \baselineskip20pt
\parskip10pt
\parindent.4in
\begin{abstract} \noi  \textcolor{red}{ Paper forthcoming in  Stochastic Analysis and Applications}\\  The finite state semi-Markov process is a generalization over the Markov chain  in which the sojourn time distribution is any general distribution. In this article we provide a sufficient stochastic maximum principle for the optimal control of a semi-Markov modulated  jump-diffusion process in which the drift, diffusion and the jump kernel of the jump-diffusion process is modulated by a semi-Markov process.  We also connect the sufficient stochastic maximum principle with the dynamic programming equation.   We apply our results to finite horizon risk-sensitive control portfolio optimization problem and to a quadratic loss minimization problem.
\end{abstract}

\newpage
\section{Introduction}
\indent The stochastic maximum principle is a stochastic version of the Pontryagin maximum principle which states that the any optimal control must satisfy a system of forward-backward stochastic differential {equations,} called the optimality system, and should maximize a functional, called the Hamiltonian. The converse indeed is true and gives the sufficient stochastic maximum principle. In this article we will derive sufficient stochastic maximum principle for a class of process called as the semi-Markov modulated jump-diffusion process. In this process the drift, the diffusion and the jump kernel term is modulated by an semi-Markov process.
\\ \indent An early investigation of stochastic maximum principle  and its application to finance has been credited to Cadenillas and Karatzas \cite{CK}. Framstadt et al. \cite{Fr} formulated the stochastic maximum principle for jump-diffusion process and applied it to a quadratic portfolio optimization problem. Their work has been partly generalized by Donnelly \cite{Do} who considered a Markov chain modulated diffusion process in which the drift and the diffusion term is modulated by a Markov chain.  Zhang et al. \cite{Zh} studied sufficient maximum principle of a process similar to that studied by Donnelly  additionally with a jump term whose kernel is also modulated by a Markov chain. It can be noted that the Markov modulated process has been quite popular with its recent applications to finance   for example   Options pricing (Deshpande and Ghosh \cite{DG}) and references therein and to portfolio optimization  refer Xhou and Yin \cite{XY}.  However application of semi-Markov modulated process to portfolio optimization in which the portfolio wealth process  is a semi-Markov modulated diffusion  are not many, see for example Ghosh and Goswami \cite{GG}.  Even so it appears that the sufficient maximum principle has not been formulated for  the case of a semi-Markov modulated diffusion process with jumps  and studied further in the context of quadratic portfolio optimization. Moreover, application of the sufficient stochastic maximum principle in the context of risk-sensitive control portfolio optimization with the portfolio wealth process following a semi-Markov modulated diffusion process has not been studied. This article aims to provide answers to these missing dots and connect them together.
For the same reasons,  alongwith providing a popular application of the sufficient stochastic maximum principle to a  quadratic loss minimization problem when the portfolio wealth process follows a semi-Markov modulated jump-diffusion, we also provide an example of risk-sensitive portfolio optimization for the diffusion part of the said dynamics. \\
\indent The article is organized as follows. In the next section we formally describe basic terminologies used in the article. In section 3 we detail the control problem that we are going to study. The sufficient maximum principle is proven in Section 4. This is followed by establishing  its connection with the dynamic programming. We  conclude the article by illustrating its applications to risk-sensitive control optimization and to a quadratic loss minimization problem.
\section{Mathematical Preliminaries}
We adopt the following notations that are valid  for the whole paper:\\
$\mathbb{R}$: the set of real numbers\\
$r,M$: any positive integer greater than 1.\\
$ {{\mathcal{X}=\{1,...,M\}}}.$\\
$\mathcal{C}^{1,2,1}([0,T] \times \mathbb{R}^{r} \times \mathcal{X} \times \mathbb{R}_{+})$: denote the family of all functions on $[0,T] \times \mathbb{R}^{r} \times \mathcal{X} \times \mathbb{R}_{+}$ which are twice continuously differentiable in $x$ and continuously differentiable in $t$ and $y$.\\
$v^{'}$, $A^{'}$: the transpose of the vector  (say )$v$ and matrix say $A$   respectively.\\
$||v||$: Euclidean norm of a vector $v$.\\
$|A|$:  norm of a matrix $A$.\\
$tr(A)$: trace of a square matrix $A$.\\
$C^{m}_{b}(\mathbb{R}^{r})$: Set of real $m$-times continuously differentiable functions which are bounded together with their derivatives upto the $m^{th}$ order.  \\
\indent We assume that the probability space ($\Omega,\mathcal{F},\{\mathcal{F}({t})\},\mathbb{P}$) is complete with filtration
$\{\mathcal{F}({t})\}_{t \geq 0}$ and is right-continuous and $\mathcal{F}({0})$ contains all $\mathbb{P}$ null sets.  Let $\{{\theta}({t})\}_{t\geq 0}$ be a semi-Markov process taking values in $ {\mathcal{X}}$ with transition probability $ {p_{ij}}$ and conditional holding time distribution $F^{h}(t|i)$.  Thus if $0 \leq t_{0}\leq t_{1}\leq ...$ are times when jumps occur, then
\begin{eqnarray}\label{2.1}
P(\theta({t_{n+1}})=j,t_{n+1}-t_{n} \leq t|\theta({t_{n}})=i)=p_{ij}F^{h}(t|i).
\end{eqnarray}
Matrix $[p_{ij}]_{\{i,j=1,...,M\}}$ is irreducible and for each $i$, $F^{h}(\cdot|i)$ has continuously differentiable and bounded density $f^{h}(\cdot|i)$. For a fixed $t$, let $n(t) \triangleq \max\{n: t_n \leq t\}$ and $Y(t) \triangleq t- t_{n(t)}$. Thus $Y(t)$ represents the amount of time the proess $\theta(t)$ is at the current state after the last jump.  The process ($\theta{(t)},Y{(t)}$)defined on ($\Omega,\mathcal{F},\mathbb{P}$)  is jointly Markov and the differential generator $\mathcal{L}$ given as follows (Chap.2, \cite{GS})
\begin{eqnarray}\label{2.3}
\mathcal{L}\phi(i,y)=\frac{d}{dy}\phi(i,y)+\frac{f^{h}(y|i)}{1-F^{h}(y|i)}\sum_{j \neq i,j \in \mathcal{X}}{p_{ij}[\phi(j,0)-\phi(i,y)]}.
\end{eqnarray}
for $\phi:\mathcal{X} \times \mathbb{R_{+}}\rightarrow \mathbb{R}$ is  { $C^{1}$} function.\\
\indent We first represent semi-Markov process $\theta(t)$ as a stochastic integral with respect to a Poisson random measure. With that perspective in mind, embed $\mathcal{X}$ in $\mathbb{R}^{M}$ by identifying $i$ with $e_{i} \in \mathbb{R}^{M}$. For $y \in [0,\infty)$ $i, j \in \mathcal{X}$, define
\begin{eqnarray*}
\lambda_{ij}(y)&=&p_{ij}\frac{f^{h}(y/i)}{1-F^{h}(y/i)} \geq 0 ~~\mbox{and}~~ \forall~~ i \neq j, \\
\lambda_{ii}(y)&=&-\sum_{j\in \mathcal{X},j \neq i}^{M}{\lambda_{ij}(y)}~~ \forall~~ i~~ \in \mathcal{X}.
\end{eqnarray*}
\indent For $i \neq j \in \mathcal{X}$ , $y \in \mathbb{R}_{+}$ let $\Lambda_{ij}(y)$ be consecutive (with respect to lexicographic ordering on $\mathcal{X}\times \mathcal{X}$) left-closed, right-open intervals of the real line, each having length $\lambda_{ij}(y)$. Define the functions $\bar{h}:\mathcal{X}\times \mathbb{R}_{+}\times\mathbb{R}\rightarrow \mathbb{R}^{r}$ and $\bar{g}:\mathcal{X}\times \mathbb{R}_{+}\times \mathbb{R} \rightarrow \mathbb{R}_{+}$ by
$$
\bar{h}(i,y,z) = \left\{ \begin{array}{rl}
 j-i &\mbox{ if $z \in \Lambda_{ij}(y)$} \\
  0 &\mbox{ otherwise}
       \end{array} \right.
$$
$$
\bar{g}(i,y,z) = \left\{ \begin{array}{rl}
 y &\mbox{ if $z \in \Lambda_{ij}(y), j \neq i$} \\
  0 &\mbox{ otherwise}
       \end{array} \right.
$$
\\
\indent Let $\mathcal{M}(\mathbb{R}_{+} \times \mathbb{R})$ be the set of all nonnegative integer-valued $\sigma$-finite measures on  { Borel } $\sigma$-field of ($\mathbb{R}_{+}\times \mathbb{R}$).
{The process $\{\tilde{\theta}{(t)},Y{(t)}\}$ is defined} by the following stochastic integral equations:
\begin{eqnarray}\label{2.2}
\begin{split}
\tilde{\theta}{(t)}=\tilde{\theta}{(0)}+\int_{0}^{t}\int_{\mathbb{R}}{\bar{h}(\tilde{\theta}{(u-)},Y{(u-)},z)N_{1}(du,dz)},\\
Y{(t)}=t-\int_{0}^{t}\int_{\mathbb{R}}{\bar{g}(\tilde{\theta}{(u-)},Y{(u-)},z)N_{1}(du,dz)},
\end{split}
\end{eqnarray}
where $N_{1}(dt,dz)$ is an $\mathcal{M}$($\mathbb{R}_{+}\times \mathbb{R}$)-valued Poisson random measure with intensity $dt m(dz)$
independent of  the $\mathcal{X}$-valued random variable  $\tilde{\theta}{(0)}$, where $m(\cdot)$ is a Lebesgue measure on $\mathbb{R}$. As usual by definition $Y(t)$ represents the amount of time, process $\tilde{\theta}(t)$ is at the current state after the last jump. We define  the corresponding compensated or centered one dimensional Poisson measure as $\tilde{N}_{1}(ds,dz)=N_{1}(ds,dz)-dsm(dz)$. It  was shown in Theorem 2.1 of Ghosh and Goswami \cite{GG} that $\tilde{\theta}{(t)}$ is a semi-Markov process with transition probability matrix $[p_{ij}]_{\{i,j=1,...,M\}}$ with conditional holding time distributions $F^{h}(y|i)$. {Since by definition $\theta(t)$ is also a semi-Markov process with transition probability matrix $[p_{ij}]_{\{i,j=1,...,M\}}$ with conditional holding time distributions $F^{h}(y|i)$ defined on the same underlying probability space, by equivalence, $\tilde{\theta}{(t)}=\theta{(t)}$ for $t \geq 0$}.\\
{
{\bf Remark 2.1}~~The semi-Markov process with conditional density $f^{h}(y|i)=\tilde{\lambda}_{i}e^{-\tilde{\lambda}_{i}y}$ for some $\tilde{\lambda}_{i}>0$, $i =1,2...,M$,  is infact a Markov chain.}
\section{The control problem}
Let $\mathcal{U} \subset \mathbb{R}^{r}$ be a closed subset. { Let $\mathbb{B}_{0}$ be the family of Borel sets  $\Gamma \subset  \mathbb{R}^{r}$ whose closure $\bar{\Gamma}$ does not contain {0}. For and Borel set $B \subset \Gamma$, one dimensional poisson random measure $ N(t,B)$ counts the number of jumps on $[0,t]$ with values in $B$.} { For a predictable process $u:[0,T] \times \Omega \rightarrow \mathcal{U}$ with right continuous left limit paths, consider the controlled process $X$ with given initial condition $X(0)=x \in \mathbb{R}^{r}$ given by}
\begin{eqnarray}\label{3.1}
dX({t})=b(t,X({t}),u({t}),\theta({t}))dt+\sigma(t,X({t}),u({t}),\theta({t}))dW({t})+\int_{\Gamma}g(t,X({t}),u({t}),\theta({t})),\gamma){N}(dt,d\gamma),\no\\
\end{eqnarray}
where $X(t) \in \mathbb{R}^{r}$ and $W(t)=(W_{1}(t),...,W_{r}(t))$ is $r$-dimensional  standard  Brownian motion. The coefficients $b(\cdot,\cdot,\cdot,\cdot):[0,T] \times \mathbb{R}^{r}\times \mathcal{U} \times \mathcal{X} \rightarrow \mathbb{R}^{r}$,$\sigma(\cdot,\cdot,\cdot,\cdot):[0,T] \times \mathbb{R}^{r}\times \mathcal{U}\times \mathcal{X} \rightarrow \mathbb{R}^{r} \times \mathbb{R}^{r}$ and $g(\cdot,\cdot,\cdot,\cdot,\cdot):[0,T] \times \mathbb{R}^{r}\times \mathcal{U} \times \mathcal{X} \times \Gamma \rightarrow \mathbb{R}^{r}$ { and satisfy the following conditions,\\
{\bf Assumption (A1)}\\
\textit{(At most linear growth)~~ There exists a constant $ C_{1}< \infty $ for any $ i \in \mathcal{X} $ such that}\\
${|\sigma(t,x,u,i)|}^{2}  +{||b(t,x,u,i)||}^{2}+\int_{\mathbb{R}}{{||g(t,x,u,i, \gamma)||}^{2}}\lambda(d\gamma) \leq C_{1}(1+||x||^{2})$\\
\textit{(Lipschitz continuity)~~ There exists a constant $C_{2}< \infty$ for any $ i \in \mathcal{X} $ such that}\\
${|\sigma(t,x,u,i)-\sigma(t,y,u,i)|}^{2} +{||b(t,x,u,i)-b(t,y,u,i)||}^{2}+\int_{\Gamma}{||g(t,x,u,i,\gamma)-g(t,y,u,i,\gamma)||^{2}}\lambda(d\gamma) \leq C_{2}||x-y||^{2}$
$\forall x,y \in \mathbb{R}^{r}$.\\
Then $X(t)$ is a unique cadlag adapted solution given by (\ref{3.1}) refer Theorem 1.19 of \cite{Oks}.}\\
\indent Define $a(t,x,u,i)=\sigma(t,x,u,i)\sigma'(t,x,u,i)$ is a $\mathbb{R}^{r\times r}$ matrix and $a_{kl}(t,x,u,i)$ is the $(k,l)^{th}$ element of the matrix $a$ while $b_{k}(t,x,u,i)$ is the $k^{th}$ element of the vector $b(t,x,u,i)$.
 We assume that
$N(\cdot,\cdot), N_{1}(\cdot,\cdot)$ and $\theta_{0},W_{t},X_{0}$ defined on ($\Omega,\mathcal{F},\mathbb{P}$) are independent. For future use we define the compensated Poisson measure $\tilde{N}(dt,d\gamma)=N(dt,d\gamma)-{\lambda} \pi(d\gamma)dt$, where $\pi(\cdot)$ is the jump distribution { (is a probability measure)  and  $0<{\lambda}<\infty$ is the jump rate} { such that $\int_{\Gamma}{\min({||\gamma||}^{2},1)}\lambda{(d\gamma)}<\infty$.}\\ \indent Consider the performance criterion
\begin{eqnarray}\label{3.2}
J^{u}(x,i,y)=E^{x,i,y}[\int_{0}^{T}{f_{1}(t,X({t}),u(t),\theta({t}),Y(t))dt+f_{2}(X(T),\theta(T),Y(T))}],
\end{eqnarray}
where $f_{1}:[0,T] \times \mathbb{R}^{r}\times \mathcal{U} \times \mathcal{X} \times \mathbb{R}_{+} \rightarrow \mathbb{R}$ is continuous and $f_{2}: \mathbb{R}^{r} \times \mathcal{X} \times \mathbb{R}_{+}\rightarrow \mathbb{R}$ is concave. We say that the admissible class of controls $u \in \mathcal{A}(T)$ if
\begin{eqnarray*}
E^{x,i,y}\bigg[\int_{0}^{T}|f_{1}(t,X(t),u(t),\theta(t),Y(t))|dt+f_{2}(X(T),\theta(T),Y(T))]\bigg]<\infty.
\end{eqnarray*}
The problem is to maximize $J^{u}$ over all $u \in \mathcal{A}(T)$ i.e. we seek $\hat{u} \in \mathcal{A}(T)$ such that
\begin{eqnarray}\label{3.3}
J^{\hat{u}}(x,i,y)=\sup_{u \in \mathcal{A}(T)}J^{u}(x,i,y),
\end{eqnarray}
where $\hat{u}$ is an optimal control.\\
Define a Hamiltonian $\mathcal{H}: [0,T] \times \mathbb{R}^{r}\times \mathcal{U} \times \mathcal{X} \times \mathbb{R}_{+} \times \mathbb{R}^{r} \times \mathbb{R}^{r \times r} \times \mathbb{R}^{r} \rightarrow \mathbb{R}$ by,
\begin{eqnarray}\label{3.4}
\mathcal{H}(t,x,u,i,y,p,q,\eta)&:=& f_{1}(t,x,u,i,y)+\bigg(b^{'}(t,x,u,i)-\int_{\Gamma}{g^{'}(t,x,u,i,\gamma)}\pi(d\gamma)\bigg)p+tr(\sigma^{'}(t,x,u,i)q)\no\\
&+&\bigg(\int_{{\Gamma}}{g^{'}(t,x,u,i,\gamma)}\pi(d\gamma)\bigg)\eta.
\end{eqnarray}
 We assume that the Hamiltonian $\mathcal{H}$ is differentiable with respect to $x$. {The adjoint equation corresponding to $u$ and $X^{u}$ in the unknown adapted processes $p(t) \in \mathbb{R}^{r}$,$ q(t) \in \mathbb{R}^{r \times r}$, $\eta:\mathbb{R}_{+} \times \mathbb{R}^{r}-\{0\}\rightarrow \mathbb{R}^{r }$ and $\tilde{\eta}(t,z)=(\eta^{(1)}(t,z),...,\eta^{(r)}(t,z))^{'}$, where $\tilde{\eta}^{(n)}(t,z) \in \mathbb{R}^{r \times r}$ for each $n=1,2,...,r$, is the backward stochastic differential equation (BSDE)},
\begin{eqnarray}\label{3.5}
dp(t)&=& -\nabla_{x}\mathcal{H}(t,X(t),u(t),\theta(t),p(t),q(t),\eta(t,\gamma))dt+q^{'}(t)dW(t)+\int_{\Gamma}{\eta(t,\gamma) \tilde{N}(dt,d\gamma)}\no \\&+&\int_{\mathbb{R}}\tilde{\eta}(t,z) \tilde{N}_{1}(dt,dz), \no \\
p(T)&=&\nabla_{x}f_{2}(X(T),\theta(T),Y(T)).~~a.s.
\end{eqnarray}
We have assumed that $\mathcal{H}$ is differentiable with respect to $x=X(t)$ and is denoted as \\$\nabla_{x}\mathcal{H}(t,X(t),u(t),\theta(t),p(t),q(t),\eta(t,\gamma))$. { As  per Remark 2.1, for the special case where the semi-Markov process has exponential holding time distribution, we would  have (\ref{3.5}) to be a BSDE with Markov chain switching. For this special case, Cohen and Elliott \cite{CE} have provided conditions for uniqueness of the solution. However, corresponding uniqueness  result  for the semi-Markov modulated BSDE as in (3.5)  seems not available in the literature. Since this paper concerns sufficient conditions, we will assume ad hoc that a solution to this BSDE  exists and is unique. }\\
{\bf Remark 3.1}~~ Notice that there are jumps in the adjoint equation (3.5)  attributed to jumps in the semi-Markov process  $\theta({t})$. This is because the drift, the diffusion and the jump kernel of the process $X({t})$ is modulated by a semi-Markov process. Also note that  the unknown process $\tilde{\eta}(t,z)$ in the adjoint equations (\ref{3.5}) does not appear in the Hamiltonian (\ref{3.4}).
\section{Sufficient Stochastic  Maximum principle}
In this section we state and prove the sufficient stochastic maximum principle.\\
{\bf Theorem 4.1}(Sufficient Maximum principle) Let $\hat{u} \in \mathcal{A}(T)$ with corresponding solution $\hat{X} \triangleq X^{\hat{u}}$. Suppose there exists a solution ($\hat{p}(t),\hat{q}(t),\hat{\eta}(t,\gamma),\hat{\tilde{\eta}}(t,z)$)of the adjoint equation (\ref{3.5}) satisfying \\
\begin{eqnarray}\label{4.1}
&&E \int_{0}^{T}{||\bigg(\sigma(t,\hat{X}(t),\theta(t))-\sigma(t,X^{u}(t),\theta(t))\bigg)^{'}\hat{p}(t)||^{2}}dt< \infty \\
&&E \int_{0}^{T}{||\hat{q}^{'}(t)\bigg(\hat{X}(t)-X^{u}(t)\bigg)||^{2}}dt< \infty \\
&&E \int_{0}^{T}{||(\hat{X}(t)-X^{u}(t))^{'}\hat{\eta}(t,\gamma)||^{2}\pi(d\gamma)}dt< \infty \\
&&E \int_{0}^{T}{|\bigg(\hat{X}(t)-X^{u}(t)\bigg)^{'}\hat{\tilde{\eta}}(t,z)|^{2}m(dz)}dt< \infty.
\end{eqnarray}
for all admissible controls $u \in \mathcal{A}(T)$. If we further suppose that \\
1.
\begin{eqnarray}\label{4.5}
\mathcal{H}(t,\hat{X}({t}),\hat{u}({t}),\theta(t),Y(t),\hat{p}({t}),\hat{q}({t}),\hat{\eta}{(t,\cdot)})=\sup_{u \in\mathcal{A}(T)} \mathcal{H}(t,\hat{X}({t}),{u}({t}),\theta(t),Y(t),\hat{p}({t}),\hat{q}({t}),\hat{\eta}{(t,\cdot)}).
\end{eqnarray}
2. for each fixed pair $(t,i,y) \in ([0,T] \times \mathcal{X} \times \mathbb{R}_{+})$,~~$\hat{\mathcal{H}}(x):= \sup_{ u \in \mathcal{A}(T)}\mathcal{H}(t,x,u,i,y,\hat{p}(t),\hat{q}(t),\hat{\eta}(t,\cdot))$ exists and is a concave function of $x$. Then $\hat{u}$ is an optimal control.\\
{\textit{Proof}}~~Fix $u \in \mathcal{A}(T)$ with corresponding solution $X=X^{u}$. For sake of brevity we would henceforth represent ($t,\hat{X}(t-),\hat{u}(t-),\theta(t-),Y(t-)$) by ($t,\hat{X}(t-)$) and ($t,{X}(t-),{u}(t-),\theta(t-),Y(t-)$) by ($t,{X}(t-)$). Then,
\begin{eqnarray*}
J(\hat{u})-J(u)=E\bigg(\int_{0}^{T}\bigg({f_{1}(t,\hat{X}(t))-f_{1}(t,X(t))}\bigg)dt+f_{2}(\hat{X}(T),\theta(T),Y(T))-f_{2}(X(T),\theta(T),Y(T))\bigg).
\end{eqnarray*}
 By use of concavity of $f_{2}(\cdot,i,y)$ we have for each $i \in \mathcal{X},~ y \in \mathbb{R}_{+}$ and (\ref{3.5}) to obtain the inequalities,
\begin{eqnarray*}
E\bigg(f_{2}(\hat{X}(T),\theta(T),Y(T))-f_{2}({X}(T),\theta(T),Y(T))\bigg) & \geq &
E\bigg((\hat{X}(T)-X(T))^{'}\nabla_{x}f_{2}(\hat{X}(T),\theta(T),Y(T))\bigg) \no \\
&\geq & E \bigg((\hat{X}(T)-X(T))^{'}\hat{p}(T)\bigg).
\end{eqnarray*}
which gives
\begin{eqnarray}\label{4.6}
J(\hat{u})-J(u) \geq E {\int_{0}^{T}{\bigg(f_{1}(t,\hat{X}(t))-f_{1}(t,X(t))\bigg)}}dt
+ E\bigg((\hat{X}(T)-X(T))^{'}\hat{p}(T)\bigg).
\end{eqnarray}
We now expand the above equation (\ref{4.6}) term by term. For the first term in this equation we use the definition of $\mathcal{H}$ as in (\ref{3.4}) to obtain
\begin{eqnarray}\label{4.7}
&&E\int_{0}^{T}{\bigg(f_{1}(t,\hat{X}(t))-f_{1}(t,X(t))\bigg)}dt \no\\
&=& E\int_{0}^{T}\bigg(\mathcal{H}(t,\hat{X}(t),\hat{u}(t),\theta(t),\hat{p}(t),\hat{q}(t),\hat{\eta}(t,\gamma))\no \\
&-&\mathcal{H}(t,{X}(t),{u}(t),\theta(t),{p}(t),{q}(t),{\eta}(t,\gamma))\bigg)dt \no \\
&-&E\int_{0}^{T}\bigg[\bigg(b(t,\hat{X}(t))-b(t,{X}(t))\no\\
&-&\int_{\Gamma}{\bigg(g(t,\hat{X}(t-),\hat{u}(t-),\theta(t-),\gamma)-g(t,{X}(t-),{u}(t-),\theta(t-),\gamma)\bigg)}\pi(d\gamma)\bigg)\hat{p}(t)\no\\
&+&tr\bigg((\sigma(t,\hat{X}(t))-\sigma(t,X(t)))^{'}\hat{q}(t)\bigg)\no \\
&+&\int_{\Gamma}(g(t,\hat{X}(t-),\hat{u}(t-),\theta(t-),\gamma)-g(t,{X}(t-),{u}(t-),\theta(t-),\gamma))^{'}\eta(t,\gamma)\pi(d\gamma)\bigg]dt. \no \\
\end{eqnarray}
To expand the second term on the right hand side of (\ref{4.6}) we begin by applying the integration by parts formula to get,
\begin{eqnarray*}
(\hat{X}(T)-X(T))^{'}\hat{p}(T)&=&  \int_{0}^{T}{(\hat{X}(t)-X(t))^{'}}d\hat{p}(t)\\
&+&\int_{0}^{T}{\hat{p}^{'}(t)d(\hat{X}(t)-X(t))}+[\hat{X}-X,\hat{p}](T).
\end{eqnarray*}
Substitute for $X$, $\hat{X}$ and $\hat{p}$ from (\ref{3.1}) and (\ref{3.5}) to obtain,
\begin{eqnarray*}
&&(\hat{X}(T)-X(T))^{'}\hat{p}(T) \\ &=&\int_{0}^{T}(\hat{X}(t)-X(t))^{'}\bigg(-\nabla_{x}\mathcal{H}(t,\hat{X}({t}),\hat{u}({t}),\hat{p}({t}),\hat{q}(t),\hat{\eta}({t,\gamma}))dt+\hat{q}^{'}(t)dW(t)\\
&+&\int_{\Gamma}{\hat{\eta}(t,\gamma)\tilde{N}(dt,d\gamma)}+\int_{\mathbb{R}}{\hat{\tilde{\eta}}(t,z)\tilde{N}_{1}(dt,dz)}\bigg)\\
&+&\int_{0}^{T}\hat{p}^{'}(t)\bigg\{\bigg(\bigg(b(t,\hat{X}(t))-b(t,X(t))\bigg)
-\int_{\Gamma}\bigg(g(t,\hat{X}(t),\hat{u}(t-),\theta({t-}),\gamma)\\
&-&g(t,{X}(t-),u(t-),\theta({t-}),\gamma)\bigg)\pi(d\gamma)\bigg)dt\\
&+&\bigg(\sigma(t,\hat{X}(t))-\sigma(t,X(t))\bigg)^{'}dW(t)\\
&+&\int_{\Gamma}{\bigg(g(t,\hat{X}(t-),\hat{u}(t-),\theta({t-}),\gamma)-g(t,{X}(t-),u(t-),\theta({t-}),\gamma)\bigg)}\tilde{N}(dt,d\gamma)\bigg\}\\
&+&\int_{0}^{T}\bigg[tr\bigg(\hat{q}^{'}(t)\bigg(\sigma(t,\hat{X}(t))-\sigma(t,X(t))\bigg)\bigg)\\
&+&\int_{\Gamma}\bigg({g(t,\hat{X}(t),\hat{u}(t-),\theta({t-}),\gamma)-g(t,{X}(t),u(t-),\theta({t-}),\gamma)}\bigg)^{'}\eta(t,\gamma)\pi(d\gamma)\bigg]dt.
\end{eqnarray*}
Due to integrability conditions ({4.1})-({4.4}), the integral with respect to the Brownian motion and the Poisson random measure are square integrable martingales which are null at the  origin. Thus taking expectations we obtain
\begin{eqnarray*}
E\bigg((\hat{X}(T)&-&X(T))^{'}\hat{p}(T)\bigg) \\ &=&\int_{0}^{T}(\hat{X}(t)-X(t))^{'}\bigg(-\nabla_{x}\mathcal{H}(t,\hat{X}({t}),\hat{u}({t}),\hat{p}({t}),\hat{q}(t),\hat{\eta}({t,\gamma}))\bigg)dt\\
&+&\int_{0}^{T}\bigg[\hat{p}^{'}(t)\bigg(b(t,\hat{X}(t))-b(t,X(t))
-\int_{\Gamma}\bigg(g(t,\hat{X}(t-),\hat{u}(t-),\theta({t-}),\gamma)\\
&-&g(t,{X}(t),u(t-),\theta({t-}),\gamma)\bigg)\pi(d\gamma)\bigg)\\
&+&\int_{0}^{T} tr\bigg(\hat{q}^{'}(t)(\sigma(t,\hat{X}(t))-\sigma(t,X(t)))\bigg)\\
&+&\int_{\Gamma}{\bigg(\bigg({g(t,\hat{X}(t-),\theta({t-}),u(t-),\gamma)-g(t,{X}(t-),\theta({t-}),u(t-),\gamma)}\bigg)^{'}\eta(t,\gamma))\bigg)\pi(d\gamma)}\bigg]dt.
\end{eqnarray*}
\begin{eqnarray*}
\end{eqnarray*}
Substitute the last equation and (\ref{4.7}) into the inequality (\ref{4.6}) to find after cancellation that
\begin{eqnarray}\label{4.8}
J(\hat{u})-J(u) & \geq &
E\int_{0}^{T}\bigg(\mathcal{H}(t,\hat{X}(t),\hat{u}(t),\theta(t),\hat{p}(t),\hat{q}(t),\hat{\eta}(t,\gamma))-\mathcal{H}(t,{X}(t),{u}(t),\theta(t),{p}(t),{q}(t),\eta(t,\gamma))\no\\
&-&(\hat{X}(t)-X(t))^{'}\nabla_{x}\mathcal{H}(t,\hat{X}(t),\hat{u}(t),\theta(t),\hat{p}(t),\hat{q}(t),\hat{\eta}(t,\gamma))\bigg)dt.
\end{eqnarray}
We can show that the integrand on the RHS of (\ref{4.8}) is non-negative a.s. for each $t \in [0,T]$ by fixing the state of the semi-Markov process and then using the assumed concavity of $\hat{\mathcal{H}}(x)$, we apply the argument in Framstad et al. \cite{Fr} . This gives $J(\hat{u}) \geq J(u)$ and $\hat{u}$ is an optimal control.$\qed$
\section{Connection to the Dynamic programming}
We show the connection between the stochastic maximum principle and dynamic programming  principle for the semi-Markov modulated regime switching jump diffusion. This tantamounts to explicitly showing connection between the value function $V(t,x,i,y)$ of the control problem and the adjoint processes $p(t), q(t)$ ,$\eta(t,\gamma)$ and $\tilde{\eta}(t,z)$. In order to apply the dynamic programming principle we put the problem into a Markovian framework by defining
\begin{eqnarray}\label{5.1}
J^{u}(t,x,i,y) \triangleq E^{X(t)=x,\theta(t)=i,Y(t)=y}[\int_{t}^{T}{f_{1}(t,X({t}),u(t),\theta({t}),Y({t}))dt+f_{2}(X(T),\theta(T),Y(T))}].
\end{eqnarray}
and put
\begin{eqnarray}\label{5.2}
V(t,x,i,y)=\sup_{u \in \mathcal{A}(T)}J^{u}(t,x,i,y)~~~~\forall~~(t,x,i,y) \in [0,T] \times \mathbb{R}^{r} \times \mathcal{X}\times \mathbb{R}_{+}.
\end{eqnarray}\\
{\bf Theorem 5.1}~~\textit{Assume that $V(\cdot,\cdot,i,\cdot)\in \mathcal{C}^{1,{3},1}([0,T]\times \mathbb{R}^{r}\times \mathcal{X}\times \mathbb{R}_{+})$ for each $i,j \in \mathcal{X}$ and that there exists an optimal Markov control $\hat{u}(t,x,i,y)$ for (\ref{5.2}), with the corresponding solution $\hat{X}=X^{(\hat{u})}$. Define
\begin{eqnarray}\label{5.3}
p_{k}(t) &\triangleq & \frac{\partial V}{\partial x_{k}}(t,\hat{X}(t),\theta(t),Y(t)).
\end{eqnarray}
\begin{eqnarray}\label{5.4}
q_{kl}(t) &\triangleq & \sum_{i=1}^{r}{\sigma_{il}(t,\hat{X}(t),\hat{u}(t),\theta(t))\frac{\partial^{2} V}{{\partial x_{i}}{\partial x_{k}}}(t,\hat{X}(t),\theta(t),Y(t))}.
\end{eqnarray}
\begin{eqnarray}\label{5.5}
\eta^{(k)}(t,\gamma) &\triangleq & \frac{\partial V}{\partial x_{k}}(t,\hat{X}(t),j,Y(t))-\frac{\partial V}{\partial x_{k}}(t,\hat{X}(t),i,Y(t)).
\end{eqnarray}
\begin{eqnarray}\label{5.6}
\tilde{\eta}^{(k)}(t,z) &\triangleq & \frac{\partial V}{\partial x_{k}}(t,\hat{X}(t-),\theta(t-)+\bar{h}(\theta({t-}),Y({t-}),z),Y({t-})-\bar{g}(\theta({t-}),Y({t-}),z))\no \\
&-&\frac{\partial V}{\partial x_{k}}(t,\hat{X}({t-}),\theta({t-}),Y({t-})).
\end{eqnarray}
{for each $(k,l =1,...,r)$. Also we assume that the coefficients $b(t,x,u,i)$, $\sigma(t,x,u,i)$ and $g(t,x,u,i,\gamma)$ belong to $C^{1}_{b}(\mathbb{R}^{r})$.} Then $p(t), q(t), \eta(t,\gamma)$ and $\tilde{\eta}(t,z)$ solves the adjoint equation (\ref{3.5}).}\\\\
We prove this theorem by using the following Ito's formula.\\
{\bf Theorem 5.2}~~{Suppose $r $ dimensional process  $X(t)=(X_{1}(t),...,X_{r}(t))$ or $\{X_{g}(t)\} $ indexed by $(g=1,2,...,r)$    satisfies the following equation,
\begin{eqnarray*}
dX_{g}(t)=b_{g}(t,X(t),u(t),\theta(t))dt+\sum_{m=1}^{r}{\sigma_{gm}(t,X(t),u(t),\theta(t))}dW_{m}(t)+\int_{\Gamma}{g_{g}(t,X(t-),u(t),\theta(t-),\gamma)}{N}(dt,d\gamma).
\end{eqnarray*}
for some $X(0)= x_{0} \in \mathbb{R}^{r}~~~a.s.$ . Further let us assume that the coefficients $b, \sigma, g$ satisfies the conditions of Assumption (A1).\\
Let $ V(\cdot,\cdot,i,\cdot)~\in~C^{1,{3},1}([0,T] \times \mathbb{R}^{r}\times \mathcal{X}\times \mathbb{R}_{+})$. Then the generalized Ito's formula is given by
\begin{eqnarray*}
&&V(t,X({t}),\theta({t}),Y({t}))- V(t,x,\theta,y)=\int_{0}^{t}{G V(s,X({s}),\theta({s}),Y({s}))ds}\\
&+&\int_{0}^{t}{(\nabla_{x} V(s,X({s}),\theta({s}),Y({s})))'\sigma(s,X({s}),\theta({s}))dW({s})}  \\
&+& \int_{0}^{t}\int_{\Gamma}[V(s,X({s-})+g(s,X({s-}),u(s),\theta({s-}),\gamma),\theta({s-}),Y({s-}))\\
&-&V(s,X({s-}),\theta({s-}),Y({s-}))]\tilde{N}(ds,d\gamma) \\
&+&\int_{0}^{t}\int_{\mathbb{R}}[V(s,X({s-}),\theta({s-})+\bar{h}(\theta({s-}),Y({s-}),z),Y({s-})-\bar{g}(\theta({s-}),Y({s-}),z))\\
 &-& V(s,X({s-}),\theta({s-}),Y({s-}))]\tilde{N}_{1}(ds,dz),
\end{eqnarray*}
where the local martingale terms are explicitly defined as \\
\begin{eqnarray*}
dM_{1}(t)&\triangleq &{(\nabla_{x} V(t,X({t}),\theta({t}),Y({t})))'\sigma(t,X({t}),u(t),\theta({t}))dW_{t}},\\
dM_{2}(t)&\triangleq &\int_{\Gamma}{[V(t,X({t-})+g(t,X({t-}),u(t),\theta({t-}),\gamma),\theta({t-}),Y({t-}))-V(t,X({t-}),\theta({t-}),Y({t-}))]\tilde{N}(dt,d\gamma)},\\
dM_{3}(t)&\triangleq &\int_{\mathbb{R}}[V\bigg(t,X({t-}),\theta({t-})+\bar{h}(\theta({t-}),Y({t-}),z),Y({t-})-\bar{g}(\theta({t-}),Y(t-),z)\bigg)\\
&-&V(t,X({t-}),\theta({t-}),Y({t-}))]\tilde{N_{1}}(dt,dz),
\end{eqnarray*}
for
\begin{eqnarray*}
G V(t,x,i,y)&=&\frac{\partial V(t,x,i,y)}{\partial t}\\&+&\frac{1}{2}\sum_{g,l=1}^{r}{a_{gl}(t,x,i)\frac{\partial V(t,x,i,y)}{\partial x_{g}\partial x_{l}}}\\&+&\sum_{g=1}^{r}{b_{g}(t,x,i)\frac{\partial V(t,x,i,y)}{\partial x_{g}}}  \\
&+&\frac{\partial V(t,x,i,y)}{\partial y}\\&+&\frac{f^{h}(y|i)}{1-F^{h}(y|i)}\sum_{j \neq i,j \in \mathcal{X},i=1}^{M}{p_{ij}[V(t,x,j,0)-V(t,x,i,y)]}  \\
&+&\lambda{\int_{\Gamma}{({V(t,x+g(t,x,i,\gamma),i,y)}-{V(t,x,i,y)})}\pi(d\gamma)},\\
\forall~t~\in~[0,T]~,x \in \mathbb{R}^{r}, (i = 1,....,M),~y \in \mathbb{R}_{+}.
\end{eqnarray*}
}
{\textit{Proof}}~~ For details refer to Theorem 5.1 in Ikeda and Watanabe \cite{IW}. $\qed$\\
{\textit{Proof of Theorem 5.1}}~~From the standard theory of the Dynamic programming the following HJB equation holds:
\begin{eqnarray*}
\frac{\partial V}{\partial t}(t,x,i,y)+\sup_{u \in \mathcal{U}}\{f_{1}(t,x,u,i,y)+\mathcal{A}^{u}V(t,x,i,y)\}=0,\\
V(T,x,i,y)=f_{2}(x,i,y).
\end{eqnarray*}
where $\mathcal{A}^{u}$ is the infinitesimal generator and the supremum is attained by $\hat{u}(t,x,i,y)$. Define
\begin{eqnarray*}
F(t,x,u,i,y)=f_{1}(t,x,u,i,y)+\frac{\partial V}{\partial t}(t,x,i,y)+\mathcal{A}^{u}V(t,x,i,y).
\end{eqnarray*}
We assume that $f_{1}$ is differentiable w.r.t to $x$. We use the Ito's formula as described in Theorem 5.2 to get,
\begin{eqnarray}
F(t,x,u,i,y)&=&f_{1}(t,x,u,i,y)+\frac{\partial V}{\partial t}(t,x,i,y)\no \\
&+&\sum_{k=1}^{r}{\frac{\partial V}{\partial x_{k}}(t,x,i,y)b_{k}(t,x,u,i) }+
\frac{1}{2}\sum_{k=1}^{r}\sum_{l=1}^{r}{\frac{\partial^{2}V}{\partial x_{k} \partial x_{l}}}(t,x,i,y)\sum_{i=1}^{r}{\sigma_{ki}(t,x,u,i)\sigma_{li}(t,x,u,i)}\no \\
&+&\sum_{j \neq i,i=1}^{M}{\frac{p_{ij}f^{h}(y|i)}{1-F^{h}(y|i)}}{(V(t,x,j,0)-V(t,x,i,y))}+\frac{\partial V}{\partial y}(t,x,i,y) \no \\
&+& \lambda \int_{\Gamma}{(V(t,x+g(t,x,u,i,\gamma),i,y)-V(t,x,i,y))}\pi(d\gamma).
\end{eqnarray}
Differentiate $F(t,x,\hat{u}(t,x,i,y),i,y)$ with respect to $x_{g}$ and evaluate at $x=\hat{X}(t)$, $i=\theta(t)$ and $y=Y(t)$, we get,
\begin{eqnarray}\label{5.8}
0&=&\frac{\partial f_{1}}{\partial x_{g}}(t,\hat{X}(t),\hat{u}(t,\hat{X}(t),\theta(t),Y(t)),\theta(t),Y(t))\no\\
&+& \frac{\partial^{2}V}{\partial x_{g} \partial t}(t,\hat{X}(t),\theta(t),Y(t))+\sum_{k=1}^{r}{\frac{\partial^{2}V}{\partial x_{g} \partial x_{k}}(t,\hat{X}(t),\theta(t),Y(t))b_{k}(t,\hat{X}(t),\hat{u}(t,\hat{X}(t),\theta(t),Y(t)),\theta(t))}\no\\
&+&\sum_{k=1}^{r}{\frac{\partial V}{\partial x_{k}}(t,\hat{X}(t),\theta(t),Y(t))\frac{\partial b_{k}}{\partial x_{g}}(t,\hat{X}(t),\hat{u}(t,\hat{X}(t),\theta(t),Y(t)),\theta(t))}\no\\
&+&\frac{1}{2}\sum_{k=1}^{r}\sum_{l=1}^{r}{\frac{\partial^{3}V}{\partial x_{g} \partial x_{k} \partial x_{l}}(t,\hat{X}(t),\theta(t),Y(t))}\no\\&\times&{\sum_{i=1}^{r}{\sigma_{k,i}(t,\hat{X}(t),\hat{u}(t,\hat{X}(t),\theta(t),Y(t)),\theta(t))\sigma_{l,i}(t,\hat{X}(t),\hat{u}(t,\hat{X}(t),\theta(t),Y(t)),\theta(t))}}\no\\
&+&\frac{1}{2}\sum_{k=1}^{r}\sum_{l=1}^{r}{\frac{\partial^{2}V}{\partial x_{k} \partial x_{l} }(t,\hat{X}(t),\hat{u}(t,\hat{X}(t),\theta(t),Y(t)),\theta(t),Y(t))}\no\\&\times&{\frac{\partial}{\partial x_{g}}\sum_{i=1}^{r}{\sigma_{k,i}(t,\hat{X}(t),\hat{u}(t,\hat{X}(t),\theta(t),Y(t)),\theta(t))\sigma_{l,i}(t,\hat{X}(t),\hat{u}(t,\hat{X}(t),\theta(t),Y(t)),\theta(t))}}\no \\
&+&\sum_{j \neq i, j \in \mathcal{X}}^{M}{\frac {p_{ij}f^{h}(y|i)}{1-F^{h}(y|i)}\bigg(\frac{\partial V}{\partial x_{g}}(t,\hat{X}(t),j,0)-\frac{\partial V}{\partial x_{g}}(t,\hat{X}(t),i,y)\bigg)}\no\\
&+&\lambda \int_{\Gamma}\bigg({\frac{\partial V}{\partial x_{g}}(t,\hat{X}(t)+g(t,\hat{X}(t),\theta(t),\gamma),\theta(t),Y(t))-\frac{\partial V}{\partial x_{g}}(t,\hat{X}(t),\theta(t),Y(t))}\bigg)\pi(d\gamma).
\end{eqnarray}
Next define, $Y_{g}=\frac{\partial V}{\partial x_{g}}(t,\hat{X}(t),\theta(t),Y(t))$ for ($g =1,...,r$). By Ito's formula (Theorem 5.2) we obtain the dynamics of $Y_{g}(t)$  as follows,
\begin{eqnarray*}
dY_{g}(t)&=&\bigg\{\frac{\partial^{2}V}{\partial {x_{g}} \partial{t}}(t,\hat{X}(t),\theta(t),Y(t))+\sum_{k=1}^{r}{\frac{\partial^{2}V}{\partial {x_{g}} \partial{x_{k}}}(t,\hat{X}(t),\theta(t),Y(t))b_{k}(t,\hat{X}(t),\hat{u}(t,\hat{X}(t),\theta(t),Y(t)),\theta(t))}\no\\
&+&\frac{1}{2}\sum_{k=1}^{r}\sum_{l=1}^{r}{\frac{\partial^{3}V}{\partial x_{g} \partial x_{k} \partial x_{l}}(t,\hat{X}(t),\theta(t),Y(t))}\no\\&\times&{\sum_{i=1}^{r}{\sigma_{ki}(t,\hat{X}(t),\hat{u}(t,\hat{X}(t),\theta(t),Y(t)),\theta(t))
\times \sigma_{li}(t,\hat{X}(t),\hat{u}(t,\hat{X}(t),\theta(t),Y(t)),\theta(t))}}\\
&+&\sum_{j \neq i, j =1}^{M}{\frac {p_{ij} f^{h}(y|i)}{1-F^{h}(y|i)}(\frac{\partial V}{\partial x_{g}}(t,\hat{X}(t),j,0)-\frac{\partial V}{\partial x_{g}}(t,\hat{X}(t),i,y))}\no \\
&+&\lambda \int_{\Gamma}{\bigg(\frac{\partial V}{\partial x_{g}}(t,\hat{X}(t)+g(t,\hat{X}(t),\hat{u}(t,\hat{X}(t),\theta(t),Y(t)),\theta(t),\gamma),\theta(t),Y(t))}\\&-&{\frac{\partial V}{\partial x_{g}}(t,\hat{X}(t),\theta(t),Y(t))\bigg)}\pi(d\gamma)\bigg\}dt\\
&+&\sum_{k=1}^{r}\frac{{\partial^{2}V}}{{\partial x_{g} \partial x_{k}}}(t,\hat{X}(t),\theta(t),Y(t))\sum_{j=1}^{r}{\sigma_{kj}(t,\hat{X}(t),\hat{u}(t,\hat{X}(t),\theta(t),Y(t)),\theta(t))dW_{j}(t)}\no \\
&+&\int_{\Gamma}\bigg\{\frac{\partial V}{\partial x_{g}}(t,\hat{X}(t-)+g(t,\hat{X}(t-),\hat{u}(t,\hat{X}(t),\theta(t),Y(t)),\theta(t-),\gamma),\theta(t-),Y(t-))\\
&-&\frac{\partial V}{\partial x_{g}}(t,\hat{X}(t-),\theta(t-),Y(t-))\bigg\}\tilde{N}(dt,d\gamma)\\
&+&\int_{\mathbb{R}}\bigg\{\frac{\partial V}{\partial x_{g}}((t,X(t-),\theta(t-)+\bar{h}(\theta(t-),Y(t-),z),Y(t-)-\bar{g}(\theta(t-),Y(t-),z)))\\&-&\frac{\partial {V}}{\partial x_{g}}(t,\hat{X}(t-),\theta(t-),Y(t-))\bigg\}{\tilde{N}_{1}(dt,dz)}.
\end{eqnarray*}
We substitute $\frac{\partial^{2}V}{\partial{x_{g}}\partial t } $ from (\ref{5.8}) to get,
\begin{eqnarray}\label{5.9}
dY_{g}(t)&=& -\frac{\partial f_{1}}{\partial x_{g}}(t,\hat{X}(t),\hat{u}(t,\hat{X}(t),\hat{u}(t,\hat{X}(t),\theta(t),Y(t)),\theta(t),Y(t)))\no\\
&-&\sum_{k=1}^{r}{\frac{\partial V}{\partial x_{k}}(t,\hat{X}(t),\theta(t),Y(t))\frac{\partial b_{k}}{\partial {x_{g}}}(t,\hat{X}(t),\hat{u}(t,\hat{X}(t),\theta(t),Y(t)),\theta(t))}\no \\
&-&\frac{1}{2}\sum_{k=1}^{r}\sum_{l=1}^{r}{\frac{\partial^{2}V}{\partial x_{k}\partial x_{l}}(t,\hat{X}(t),\hat{u}(t,\hat{X}(t),\theta(t),Y(t)),\theta(t),Y(t))}\no\\&\times&{\frac{\partial}{\partial x_{g}}(\sum_{k=1}^{r}{\sigma_{ki}(t,\hat{X}(t),\theta(t))}{ \sigma_{li}(t,\hat{X}(t),\hat{u}(t,\hat{X}(t),\theta(t),Y(t)),\theta(t))})}\no\\
&+& \sum_{k=1}^{r}{\frac{\partial^{2}V}{\partial x_{g}\partial x_{k}}(t,\hat{X}(t),\theta(t),Y(t))\sum_{j=1}^{r}{\sigma_{kj}(t,\hat{X}(t),\hat{u}(t,\hat{X}(t),\theta(t),Y(t)),\theta(t))dW_{j}(t)}}\no\\
&+&\int_{\Gamma}\bigg\{(\frac{\partial V}{\partial x_{g}}(t,\hat{X}(t-)+g(t,X(t-),\hat{u}(t,\hat{X}(t),\theta(t),Y(t)),\theta(t-),\gamma),\theta(t-),Y(t-))\no\\
&-&\frac{\partial V}{\partial x_{g}}(t,\hat{X}(t-),\theta(t-),Y(t-)))\bigg\}\tilde{N}(dt,d\gamma)\no\\
&+&\int_{\mathbb{R}}\bigg\{\frac{\partial V}{\partial x_{g}}((t,X(t-),\theta(t-)+\bar{h}(\theta(t-),Y(t-),z),Y(t-)-\bar{g}(\theta(t-),Y(t-),z)))\no\\&-&\frac{\partial {V}}{\partial x_{g}}(t,\hat{X}(t-),\theta(t-),Y(t-))\bigg\}{\tilde{N}_{1}(dt,dz)}.
\end{eqnarray}
We have the following identity,
\begin{eqnarray}\label{5.10}
&&\frac{1}{2}\sum_{k=1}^{r}\sum_{l=1}^{r}{\frac{\partial^{2}V}{\partial x_{k}\partial x_{l}}(t,\hat{X}(t),\theta(t),Y(t))}\no\\&\times&{\frac{\partial}{\partial x_{g}}\bigg(\sum_{i=1}^{r}{\sigma_{ki}(t,\hat{X}(t),\hat{u}(t,\hat{X}(t),\theta(t),Y(t)),\theta(t))\sigma_{li}(t,\hat{X}(t),\hat{u}(t,\hat{X}(t),\theta(t),Y(t)),\theta(t))}\bigg)}\no \\
&=&\sum_{k=1}^{r}\sum_{l=1}^{r}\sum_{i=1}^{r}{\sigma_{il}(t,\hat{X}(t),\hat{u}(t,\hat{X}(t),\theta(t),Y(t)),\theta(t)){\frac{\partial^{2}V}{\partial x_{i}\partial x_{k}}(t,\hat{X}(t),\theta(t),Y(t))}}\no\\&\times&{{\frac{\partial \sigma_{kl}}{\partial x_{g}}(t,\hat{X}(t),\hat{u}(t,\hat{X}(t),\theta(t),Y(t)),\theta(t))}}.
\end{eqnarray}
Next, from (\ref{3.4}) we obtain,
\begin{eqnarray}\label{5.11}
&&\frac{\partial \mathcal{H}}{\partial x_{g}}(t,X(t),u(t),\theta(t),Y(t),p(t),q(t),\eta(t,\gamma))\no \\&=&\frac{\partial f_{1}}{\partial x_{g}}(t,\hat{X}(t),\hat{u}(t,\hat{X}(t),\theta(t),Y(t)),\theta(t),Y(t))\no\\
&+&\sum_{i=1}^{r}\bigg(\frac{\partial b_{i}}{\partial x_{g}}(t,\hat{X}(t-),\hat{u}(t,\hat{X}(t),\theta(t),Y(t)),\theta(t-))\no\\&-&\int_{\Gamma}{\frac{\partial g_{i}}{\partial x_{g}}(t,X(t-),\hat{u}(t,\hat{X}(t),\theta(t),Y(t)),\theta(t-),\gamma)}\pi(d\gamma)\bigg)p_{i}(t)+tr(\frac{\partial \sigma^{'}(t,x,\hat{u},\theta(t))}{\partial x_{g}}q)\no\\
&+&{\sum_{i=1}^{r}\int_{\Gamma}{\frac{\partial g_{i}}{\partial x_{g}}(t,X(t-),\theta(t-),\gamma)}\pi(d\gamma)(\eta^{(g)}_{i}(t,\gamma))}.
\end{eqnarray}
We also note that
\begin{eqnarray*}
tr(\frac{\partial \sigma^{'}(t,x,u,i)}{\partial x_{g}}q)&=&\sum_{l=1}^{r}{[\frac{\partial \sigma^{'}(t,x,u,i)}{\partial x_{g}}q ]_{ll}}\\
&=&\sum_{l=1}^{r}\sum_{k=1}^{r}q_{k,l}\frac{\partial \sigma_{kl}}{\partial x_{g} }(t,x,u,i).
\end{eqnarray*}
Substitute (\ref{5.3})-(\ref{5.6}) and (\ref{5.11}) gives,
\begin{eqnarray}
dY_{g}(t)&=&-\frac{\partial\mathcal{H}}{\partial x_{g}}(t,X(t),u(t),\theta(t),Y(t),p(t),q(t),\eta(t,\gamma))dt+\sum_{j=1}^{r}{q_{gj}}(t)dW_{j}(t)\no\\
&+&\int_{\Gamma}\eta(t,\gamma) \tilde{N}(dt,d\gamma)+\int_{\mathbb{R}}\tilde{\eta}(t,z) \tilde{N}_{1}(dt,dz).
\end{eqnarray}
Since $Y_{g}(t)=p_{g}(t)$ for each $g=1,...,r$, we have shown that $p(t),q(t),\eta(t,\gamma)$ and $\tilde{\eta}(t,z)$  solve the adjoint equation (\ref{3.5}). $\qed$\\
\section{Applications}
We illustrate the theory developed towards applying it to some key financial wealth optimization problems. For an early motivation on applying sufficient maximum principle, we first consider wealth dynamics to follow   semi-Markov modulated diffusion (no jumps case) and apply it towards the risk-sensitive control portfolio optimization problem. We follow it up by illustrating an application of semi-Markov modulated jump-diffusion wealth dynamics to a quadratic loss minimization problem. Unless otherwise stated, all the processes defined in this section are one dimensional.\\
{\bf Risk-sensitive control portfolio optimization}~~Let us consider a financial market consisting of two continuously traded securities namely the risk less bond  and a stock. The dynamics of the riskless bond is known to follow
\begin{eqnarray*}
dS_{0}(t)=r(t,\theta(t-))S_{0}(t)dt~~~S_{0}(0)=1.
\end{eqnarray*}
where $r(t,\theta(t))$ is the risk-free interest rate at time $t$ and is modulated by an underlying semi-Markov process as described earlier.
The dynamics of the stock price is given as
\begin{eqnarray*}
dS_{1}(t)=S_{1}(t)[(\mu(t,\theta(t-)))dt+\sigma(t,\theta(t-))dW(t)],
\end{eqnarray*}
where $(\mu(t,\theta(t-)))$ is the instantaneous expected rate of return and as usual $\sigma(t,\theta(t-))$ is the instantaneous volatility rate. The stock price process is thus driven by a 1-d Brownian motion. We denote the wealth of the investor to be $X(t) \in \mathbb{R}$ at time $t$. He holds $\theta_{1}(t)$ units  of stock and $\theta_{0}(t)=1-\theta_{1}(t)$ units is held in the riskless bond market. From the self-financing principle (refer Karatzas and Shreve \cite{KS}), the wealth process follows the dynamics given as,
\begin{eqnarray*}
dX(t)=(r(t,\theta(t-))X(t)+h(t)\sigma(t,\theta(t-))\bar{m}(t,\theta(t-)))dt+h(t)\sigma(t,\theta(t-))dW(t)~~~X(0)=x,
\end{eqnarray*}
where $h(t)=\theta_{1}(t)S_{1}(t)$, { $\bar{m}(t,i)=\frac{\mu(t,i)-r(t,i)}{\sigma(t,i)} \geq 0$ and  the variables $r(t,i), b(t,i)$ and $\sigma(t,i),$ and $\sigma^{-1}(t,i)$ for each $i \in \mathcal{X}$ are measurable and uniformly bounded in $t \in [0,T]$}. { Also $h(\cdot)$ occuring in the drift and diffusion term in above dynamics of $X(t)$ satisfies the following conditions\\
1. $E[\int_{0}^{T}{h^{2}(t)dt}]< \infty$\\
2. $E[\int_{0}^{T}{|r(t,\theta(t-))X(t)+h(t)\sigma(t,\theta(t-))\bar{m}(t,\theta(t-))|}dt+\int_{0}^{T}{h^{2}(t)\sigma^{2}(t,\theta(t-))}dt]< \infty$\\
3. The SDE for $X$ has a unique strong solution.\\ These conditions on $h(\cdot)$ are needed in order to prevent doubling strategies which otherwise  would yield arbitrary  profit at time $T$ for an investor.}\\
\indent In a classical risk-sensitive control optimization problem, the investor aims to maximize over some admissible class of portfolio $\mathcal{A}(T)$ the following risk-sensitive criterion given by
\begin{eqnarray*}
J(\hat{h}(\cdot),x)&=&\max_{h \in \mathcal{A}(T)}{\frac{1}{\gamma}}\mathbb{E}[{X(T)}^{\gamma}|X(0)=x,\theta(0)=i,Y(0)=y],~~~\gamma \in (1,\infty) \\
&=& -\min_{h \in \mathcal{A}(T)} \frac{1}{\gamma} \mathbb{E}[{X(T)}^{\gamma}|X(0)=x,\theta(0)=i,Y(0)=y],
\end{eqnarray*}
where the exogenous parameter  $\gamma$ is the usual risk-sensitive criterion that describes the risk attitude of an investor. Thus the optimal expected utility function depends on $\gamma$ and is a generalization of the traditional  stochastic control approach to utility optimization in the sense that now the degree of risk aversion  of the investor is explicitly parameterized through $\gamma$ rather than importing it in the problem via an exogeneous utility function. See Whittle \cite{Wh} for general overview on risk-sensitive control optimization. We now use the sufficient maximum principle  (Theorem 4.1). Set the control problem
$u(t)\triangleq h(t)$.\\
The corresponding Hamiltonian (for the non-jump case)(\ref{3.4}) becomes,
\begin{eqnarray*}
\mathcal{H}(t,x,u,i,p,q)=(r(t,i)x+u\sigma(t,i)\bar{m}(t,i))p+u\sigma(t,i)q.
\end{eqnarray*}
The adjoint process (\ref{3.5}) is given by
\begin{eqnarray}\label{6.1}
dp(t)&=&-r(t,\theta(t-))p(t)dt+q(t)dW(t)+\int_{\mathbb{R}}{\tilde{\eta}(t,z)\tilde{N}_{1}(dt,dz)},\no \\
p(T)&=&X(T)^{\gamma-1}~~~~a.s..
\end{eqnarray}
We need to determine $p(t),q(t)$ and $\eta(t,z)$ in (\ref{6.1}). Going by the terminal condition $p(T)$ we observe that the adjoint process $p$ is the first derivative of $(x^{\gamma})$. Hence we assume that $p(t)$  defined as,
\begin{eqnarray*}
p(t)=(X(t))^{\gamma-1}e^{{\phi(t,\theta(t),Y(t))}}.
\end{eqnarray*}
where $\phi(T,\theta(T)=i,Y(T))=0~~~a.s.$ for each $i \in \{1,...,M\}$. Using the Ito's formula we get,
\begin{eqnarray}\label{6.2}
&&\frac{dp(t)}{p(t)}=\sum_{i=1}^M {1}_{\theta(t-)=i}\bigg((\gamma-1)\bigg\{(r(t,\theta(t-))+\frac{u(t)\sigma(t,\theta(t-))\bar{m}(t,\theta(t-))}{X(t)}\bigg)\no\\
&+&\frac{1}{2}(\gamma-1)(\gamma-2)\sigma^{2}(t,\theta(t-))\frac{u^{2}(t)}{X^{2}(t)} \no \\
&+&\phi_{t}(t,\theta(t-),y)+\phi_{y}(t,\theta(t-),y)+\frac{f^{h}(y|\theta(t-)=i)}{1-F^{h}(y|\theta(t-)=i)}\sum_{j \neq i}{p_{ij}(\phi(t,j,0)-\phi(t,\theta(t-),y))}\bigg\}dt \no \\
&+&{(\gamma-1)\frac{u(t)}{X(t)}\sigma(t,\theta(t-))}dW(t) \no \\ &+&\int_{\mathbb{R}}\bigg(\phi(t,X(t-),\theta(t-)+\bar{h}(\theta(t-),Y(t-),z),Y(t-)-\bar{g}(\theta(t-),Y(t-),z))\no\\
&-&\phi(t,\theta(t-),Y(t-))\bigg)\tilde{N}_{1}(dt,dz).\no\\
\end{eqnarray}
Comparing the coefficient of (\ref{6.2}) with that in (\ref{6.1}) we get
\begin{eqnarray}\label{6.3}
-r(t,\theta(t-))&=&\sum_{i=1}^M {1}_{\theta(t-)=i}\bigg((\gamma-1)\bigg(r(t,\theta(t-))+\frac{u(t)\sigma(t,\theta(t-))\bar{m}(t,i)}{X(t)}\bigg)+\frac{1}{2}(\gamma-1)(\gamma-2)\frac{u^{2}(t)}{X^{2}(t)} \no \\
&+& \phi_{t}(t,\theta(t-),y)+\phi_{y}(t,\theta(t-),y)+\frac{f^{h}(y|i)}{1-F^{h}(y|\theta(t-)=i)}\sum_{j \neq i}{p_{ij}(\phi(t,j,0)-\phi(t,\theta(t-),y))}\bigg).\no\\
\end{eqnarray}
\begin{eqnarray}\label{6.4}
{q}(t)=(\gamma-1)\frac{u(t)}{X(t)}\sigma(t,\theta(t-)){p}(t).
\end{eqnarray}
\begin{eqnarray}\label{6.5}
\tilde{\eta}(t,z)&=&\bigg(\phi(t,\theta(t-)+\bar{h}(\theta(t-),Y(t-),z),Y(t-)-\bar{g}(\theta(t-),Y(t-),z))\no \\
&-&\phi(t,\theta(t-),Y(t-))\bigg)p(t).
\end{eqnarray}
Let $\hat{u} \in \mathcal{A}(T)$ be a candidate optimal control corresponding to the wealth process $\hat{X}$ and the adjoint triplet ($\hat{p},\hat{q},\hat{\eta}$),
then from the Hamiltonian (\ref{3.4}) for all $u \in \mathbb{R}$  we have
\begin{eqnarray}\label{6.6}
\mathcal{H}(t,\hat{X}(t),u,\theta(t),\hat{p}(t),\hat{q}(t))=\bigg(r(t,\theta(t))\hat{X}(t)+u\sigma(t,\theta(t)) \bar{m}(t,\theta(t))\bigg)\hat{p}(t)+u\sigma(t,\theta(t))\hat{q}(t).
\end{eqnarray}
As this is a linear function of $u$, we guess that the coefficient of $u$ vanishes at optimality, which results in the equality
\begin{eqnarray}\label{6.7}
\bar{m}(t,\theta(t-))\hat{p}(t)+\hat{q}(t)=0.
\end{eqnarray}
Substitute  equation (\ref{6.7}) in (\ref{6.4}) to obtain the expression for the  control as
\begin{eqnarray}\label{6.8}
\hat{u}(t)=\frac{\bar{m}(t,\theta(t-))}{(1-\gamma)\sigma(t,\theta(t-))}\hat{X}(t).
\end{eqnarray}
We now aim to determine the explicit expression for ${p}(t)$ which is only possible if we can determine what $\phi(t,\theta(t),Y(t))$ is. We substitute  $\hat{u}$  from above and input it in  equation (\ref{6.3}) to get
\begin{eqnarray}\label{6.9}
0&=&\gamma r(t,\theta(t-))-{\bar{m}^{2}(t,\theta(t-))}+\frac{(2-\gamma)}{(1-\gamma)}\frac{\bar{m}^{2}(t,\theta(t-))}{2\sigma^{2}(t,\theta(t-))}\no\\
&+& \phi_{t}(t,\theta(t-),y)+\phi_{y}(t,\theta(t-),y)+\frac{f^{h}(y|\theta(t-)=i)}{1-F^{h}(y|\theta(t-)=i)}\sum_{i=1,j \neq i}^{M}{p_{ij}(\phi(t,j,0)-\phi(t,\theta(t-),y))}.\no\\
\end{eqnarray}
with terminal boundary condition given as $\phi(T,\theta(T),Y(T))=0$~~~a.s. Consider the process
\begin{eqnarray}\label{6.10}
\tilde{\phi}(t,\theta(t),Y(t)) \triangleq
 E\bigg[\exp\bigg(\int_{t}^{T}\bigg\{\gamma r(s,\theta(s))-{\bar{m}^{2}(s,\theta(s))}+\frac{(2-\gamma)}{(1-\gamma)}\frac{\bar{m}^{2}(s,\theta(s))}{{2\sigma^{2}(s,\theta(s))}}\bigg\}ds\bigg)|{\theta(t-)=i,Y(t-)=y}\bigg].\no\\
\end{eqnarray}
We aim to show that $\phi=\tilde{\phi}$. For the same we define the following martingale,
\begin{eqnarray}\label{6.11}
R(t)\triangleq E\bigg[\exp\bigg(\int_{0}^{T}{\bigg\{\gamma r(s,\theta(s))-{\bar{m}^{2}(s,\theta(s))}+\frac{(2-\gamma)}{(1-\gamma)}\frac{\bar{m}^{2}(s,\theta(s))}{{2\sigma^{2}(s,\theta(s))}}\bigg\}}ds\bigg)|\mathcal{F}_{t}^{\theta,y}\bigg],
\end{eqnarray}
where $\mathcal{F}_{\tau}^{\theta,y} \triangleq \sigma\{\theta{(\tau)},Y(\tau), \tau \in [0,t]\}$  augmented with $\mathbb{P}$ null sets is the filtration generated by the processes $\theta(t)$ and $Y(t)$. From the $\{\mathcal{F}_{t}^{\theta,y}\}$-martingale representation theorem, there exist $\{\mathcal{F}_{t}^{\theta,y}\}$-previsible, square integrable process $\nu(t,i,y)$  such that
\begin{eqnarray}\label{6.12}
R(t)=R(0)+\int_{0}^{t}\int_{\mathbb{R}}{\nu(\tau,\theta(\tau-),Y(\tau-))}\tilde{N}_{1}(d\tau,dz).
\end{eqnarray}
By positivity of $R(t)$ we can define $\hat{\nu}(\tau,\theta(\tau-),Y(\tau-)) \triangleq  (\nu(\tau,\theta(\tau-),Y(\tau-)))R^{-1}(\tau-)$ so that
\begin{eqnarray}\label{6.13}
R(t)=R(0)+{\int_{0}^{t}\int_{\mathbb{R}}{R(\tau-)\hat{\nu}(\tau,\theta(\tau-),Y(\tau-))}\tilde{N}_{1}(d\tau,dz)}.
\end{eqnarray}
From the definition of $\tilde{\phi}$ in (\ref{6.10}) and the definition of $R$ in (\ref{6.11}) it is easy to see that we have the following relationship
\begin{eqnarray}\label{6.14}
R(t)=\tilde{\phi}(t,\theta(t),Y(t))\exp\bigg\{\int_{0}^{t}(\gamma r(s,\theta(s))-{\bar{m}^{2}(s,\theta(s))}+\frac{(2-\gamma)}{(1-\gamma)}\frac{\bar{m}^{2}(s,\theta(s))}{{2\sigma^{2}(s,\theta(s))}})ds\bigg\},\no\\ ~~\forall~t~\in~[0,T].
\end{eqnarray}
Using the Ito's expansion of $\tilde{\phi}(t,\theta(t),Y(t))$  to the RHS of (\ref{6.14}) followed up by comparing it with martingale representation of $R(t)$ in (\ref{6.12}) we get $\phi:= \tilde{\phi}$. We can thus substitute $ \hat{q}$ and $\hat{\tilde{\eta}}$ in expression (\ref{6.4}),(\ref{6.5}) in lieu of $q$ and $\tilde{\eta}(t,z)$ respectively. With the choice of control $\hat{u}$ given by (\ref{6.8}) and boundedness condition on the market parameters $r,\mu$ and $\sigma$, the conditions in Theorem 4.1 are satisfied and hence $\hat{u}(t)$ is an optimal control process and the explicit representation of $\hat{p}$ is given by
\begin{eqnarray*}
\hat{p}(t)=(X(t))^{\gamma-1}e^{E[\exp(\int_{t}^{T}{\gamma r(s,\theta(s))-{\bar{m}^{2}(s,\theta(s))}+\frac{(2-\gamma)}{(1-\gamma)}\frac{\bar{m}^{2}(s,\theta(s))}{{2\sigma^{2}(s,\theta(s))}}ds|}{\theta(t-)=i,Y(t-)=y})]}.
\end{eqnarray*}
{\bf Quadratic loss minimization}~~ We now provide an example related
to quadratic loss minimization  where the portfolio wealth process is given by
\begin{eqnarray}\label{6.15}
dX^{h}({t})&=&\bigg(r({t},\theta(t))X^{h}(t)+h(t)\sigma(t,\theta(t))\bar{m}(t,\theta(t))-h(t)\int_{\Gamma}{g(t,X^{h}(t),\theta(t),\gamma)\pi(d\gamma)}\bigg)dt\no\\
&+&h(t)\sigma(t,\theta(t))dW(t)
+ h(t)\int_{\Gamma}{g(t,X^{h}(t),\theta(t),\gamma)\tilde{N}(dt,d\gamma)}, \no \\
X^{h}(0)&=&x_{0}~~a.s.
\end{eqnarray}
{where the market price of risk is defined as $\bar{m}({t},i,y)=\sigma^{-1}(t,i)(b(t,i)-r(t,i))$. {As like earlier example , we have that $\bar{m}(t,i) \geq 0$ and  that the variables $r(t,i), b(t,i)$, $\sigma(t,i)$ , $\sigma^{-1}(t,i)$ and $g(t,x,i,\gamma)$ for each $i \in \mathcal{X}$ are measurable and uniformly bounded in $t \in [0,T]$. We assume that  $g(t,x,i,\gamma)>-1$ for each $i \in \mathcal{X}$ and for a.a. $t,x,\gamma$. This insures that  $X^{h}(t)>0$ for each $t$. We further assume following conditions for each $i \in \mathcal{X}$\\
1.
$E[\int_{0}^{T}{h^{2}(t)dt}]< \infty.$\\
2.
$E[\int_{0}^{T}{|r(t,i)X(t)+h(t)\sigma(t,i)\bar{m}(t,i)|}dt+\int_{0}^{T}{h^{2}(t)\sigma^{2}(t,i)}dt+\int_{0}^{T}{h^{2}(t)g^{2}(t,X(t),i,\gamma)}dt]< \infty.$\\
3. $t \rightarrow \int_{\mathbb{R}}{h^{2}(t)g^{2}(t,x,i,\gamma)\pi(d\gamma)}$ is bounded. \\
4. the SDE for $X$ has a unique strong solution.\\
}
The portfolio process $h(\cdot)$ satisfying the above four conditions is said to be admissible and  belongs to $\mathcal{A}(T)$ (say).
} We consider the problem of finding an admissible portfolio process $h \in \mathcal{A}(T)$ such that
\begin{eqnarray*}
\inf_{h \in \mathcal{A}(T)}{E[(X^{h}(T)-d)^{2}]},
\end{eqnarray*}
over all $h \in \mathcal{A}(T)$. Set the control process $u(t) \triangleq h(t)$ and $X(t) \triangleq X^{h}(t)$.
For this example the Hamiltonian (\ref{3.4}) becomes
\begin{eqnarray}\label{6.16}
\mathcal{H}(t,x,h,i,y,p,q,\eta)&=&\bigg[r(t,i)x+u\sigma(t,i)\bar{m}(t,i)-u\int_{\Gamma}{g(t,x,i,\gamma)\pi(d\gamma)}\bigg]p+u\sigma(t,i)q \no \\ &+&\bigg(u\int_{\Gamma}{g(t,x,i,\gamma)}\pi(d\gamma) \bigg)\eta ,
\end{eqnarray}
and the adjoint equations are for all time $t \in [0,T)$,
\begin{eqnarray}\label{6.17}
dp(t)&=&-r(t,\theta(t-))p(t)dt+q(t)dW(t)+\int_{\Gamma}{\eta(t,\gamma)\tilde{N}(dt,d\gamma)}+\int_{\mathbb{R}}{\tilde{\eta}(t,z)\tilde{N}_{1}(dt,dz)}, \no \\
p(T)&=&-2X(T)+2d ~~a.s.
\end{eqnarray}
We seek to determine $p(t),q(t), \eta(t,\gamma)$ and $\tilde{\eta}(t,z)$ in (\ref{6.17}). Going by  (\ref{6.17}) we assume that ,
\begin{eqnarray}\label{6.18}
p(t)=\phi(t,\theta(t),Y(t))X(t)+\psi(t,\theta(t),Y(t)).
\end{eqnarray}
with the terminal boundary conditions being
\begin{eqnarray}\label{6.19}
\phi(T,i,y)=-2~~~~~~~~~~~~\psi(T,i,y)=2d~~~~\forall ~i~\in~\mathcal{X}.
\end{eqnarray}
For the sake of convenience we again rewrite the following Ito's formula for a function $f(t,\theta(t),y(t))\in \mathcal{C}^{1,2,1}$ given as
\begin{eqnarray}\label{6.20}
&&df(t,\theta(t),Y(t))=\bigg(\frac{\partial f(t,\theta(t),Y(t))}{\partial {t}}+\frac{(f^{h}(y/i))}{(1-F^{h}(y/i))}\sum_{j \neq i, j=1}^{M}p_{\theta(t-)=i,j}[f(t,j,0)-f(t,\theta(t-),y)]\no\\
&+&\frac{\partial f(t,\theta(t),Y(t))}{\partial y}\bigg)dt\no \\ &+&\int_{\mathbb{R}}{[f(t,\theta({t-})+\bar{h}(\theta({t-}),Y({t-}),z),Y({t-})-\bar{g}(\theta({t-}),Y({t-}),z))-f(t,\theta({t-}),Y({t-}))]\tilde{N}_{1}(dt,dz)}.\no \\
\end{eqnarray}
We apply the Ito's product rule to (\ref{6.18})  to obtain
\begin{eqnarray}\label{6.21}
dp({t})&=&X({t-})d\phi(t,\theta(t-),Y(t))+\phi(t,\theta(t-),Y(t))dX(t)+d\phi(t,\theta(t-),Y(t))dX(t)+d\psi(t)\no\\
&=& \sum_{i=1}^{M}{1_{\theta_{t-}=i}}\bigg\{X(t-)\bigg(\phi(t,\theta(t-),y)r(t,\theta(t-))+\phi_{t}(t,\theta(t-),Y(t))+\phi_{y}(t,\theta(t-),Y(t))\no \\
&+&\sum_{i=1,j \neq i}^{M}{p_{ij}\frac{f^{h}(y/i)}{1-F^{h}(y/i)}(\phi(t,j,0)-\phi(t,\theta(t-),Y(t)))}\bigg)+u(t)\phi(t,\theta(t-),Y(t))\sigma(t,\theta(t-))\bar{m}(t,\theta(t-))\no \\
&-&u(t)\phi(t,\theta(t-),Y(t))\int_{{\Gamma}}{g(t,X(t),\theta(t-),\gamma)\pi(d\gamma)}+\psi_{t}(t,\theta(t-),Y(t))+\psi_{y}(t,\theta(t-),Y(t))\no\\
&+&\sum_{i=1,i \neq j}^{M}{p_{ij}\frac{f^{h}(y/i)}{1-F^{h}(y/i)}[\psi(t,j,0)-\psi(t,\theta(t-)=i,Y(t))]}\bigg\}dt \no \\
&+&u(t)\phi(t,\theta({t-}),Y({t}))\sigma(t,\theta({t-}))dW({t})+u(t)\phi(t,\theta({t-}),Y({t-}))\int_{{\Gamma}}{g(t,X(t-),\theta({t-}),\gamma) \tilde{N}(dt,d\gamma)}\no\\
&+&\int_{\mathbb{R}}\bigg[X(t-)(\phi(t,\theta({t-})+\bar{h}(\theta({t-}),Y({t-}),z),Y({t-})-\bar{g}(\theta({t-}),Y({t-}),z))-\phi(t,\theta({t-}),Y({t-})))\no \\ &+&\psi(t,\theta({t-})+\bar{h}(\theta({t-}),Y({t-}),z),Y({t-})-\bar{g}(\theta({t-}),Y({t-}),z))-\psi(t,\theta({t-}),Y({t-}))\bigg]\tilde{N}_{1}(dt,dz). \no \\
\end{eqnarray}
Comparing coefficients with (\ref{6.17}) we obtain three equations given as
\begin{eqnarray}\label{6.22}
&-&r(t,\theta({t-}))p(t-)\no\\&=&\sum_{i=1}^{M}{1_\{{\theta_{t-}=i},Y(t-)=y}\}\bigg\{X(t-)\bigg(\phi(t,\theta({t-}),Y(t))r(t,\theta({t-}))+\phi_{t}(t,\theta({t-}),Y(t))+\phi_{y}(t,\theta({t-}),Y(t))\no\\
&+&\sum_{i=1,j \neq i}^{M}{p_{ij}\frac{f^{h}(y/i)}{1-F^{h}(y/i)}(\phi(t,j,0)-\phi(t,\theta({t-}),Y(t)))}\bigg) +u(t)\phi(t,\theta({t-}),Y(t))\sigma(t,\theta({t-}))\bar{m}(t,\theta({t-}))\no\\
&-&u(t)\phi(t,\theta({t-}),Y(t))\int_{\Gamma}{g(t,x,\theta({t-}),\gamma)\pi(d\gamma)}
+\psi_{t}(t,\theta({t-})),Y(t)+\psi_{y}(t,\theta({t-}),Y(t))\no\\
&+&\sum_{i \neq j}^{M}{p_{ij}\frac{f^{h}(y/i)}{1-F^{h}(y/i)}[\psi(t,j,0)-\psi(t,\theta({t-}),Y(t))]}\bigg\}.\no\\
\end{eqnarray}
\begin{eqnarray}\label{6.23}
q(t)=u(t)\phi(t,\theta({t-}),Y({t-}))\sigma(t,\theta({t-})).
\end{eqnarray}
\begin{eqnarray}\label{6.24}
\eta(t,\gamma)=u(t)\phi(t,\theta({t-}),Y(t-))g(t,X(t-),\theta({t-}),\gamma).
\end{eqnarray}
\begin{eqnarray}\label{6.25}
\tilde{\eta}(t,z)&=&X(t-)(\phi(t,\theta({t-})+\bar{h}(\theta({t-}),Y({t-}),z),Y({t-})-\bar{g}(\theta({t-}),Y({t-}),z))-\phi(t,\theta({t-}),Y({t-})))\no \\ &+&\psi(t,\theta({t-})+\bar{h}(\theta({t-}),Y({t-}),z),Y({t-})-\bar{g}(\theta({t-}),Y({t-}),z))-\psi(t,\theta({t-}),Y({t-})).\no\\
\end{eqnarray}
Let $\hat{u} \in \mathcal{A}(T)$ be a candidate optimal control corresponding to the wealth process $\hat{X}(T)$ and the adjoint triplet ($\hat{p},\hat{q},\hat{\eta},\hat{\tilde{\eta}}$). Then from the Hamiltonian (\ref{3.4}) for all $u \in \mathcal{A}(T)$  we have
\begin{eqnarray}\label{6.26}
\mathcal{H}(t,\hat{X}(t),u,\theta(t),\hat{p}(t),\hat{q}(t),\hat{\eta}(t))&=&\bigg(r(t,\theta(t))\hat{X}(t)+u \sigma(t,\theta(t))\bar{m}(t,\theta(t))\no\\
&-&u\int_{\Gamma}{g(t,\hat{X}(t-),\theta(t-),\gamma)\pi{d(\gamma)}}\bigg)\hat{p}(t)\no\\
&+&u\sigma(t,\theta(t))\hat{q}(t)+\bigg(u\int_{\Gamma}{g(t,\hat{X}(t-),\theta(t-),\gamma)\pi(d{\gamma})}\bigg)\hat{\eta}(t,\gamma).\no\\
\end{eqnarray}
As this is a linear function of $u$, we guess that the coefficient of $u$ vanishes at optimality, which results in the following equality
\begin{eqnarray}\label{6.27}
\hat{q}(t)&=&\bigg(-\bar{m}(t,\theta({t-}))+\frac{1}{\sigma(t,\theta({t-}))}\int_{\Gamma}{g(t,\hat{X}(t),\theta(t),\gamma)\pi(d\gamma)}\bigg)\hat{p}(t)\no\\
&-&\frac{1}{\sigma(t,\theta({t-}))}\int_{\Gamma}{(g^{'}(t,\hat{X}(t),\theta(t),\gamma))\pi(d\gamma)\hat{\eta}(t,\gamma)}.\no\\
\end{eqnarray}
Also substituting (\ref{6.27}) for $\hat{q}(t)$ in (\ref{6.23}) and using (\ref{6.18}) and(\ref{6.24}) we get,
\begin{eqnarray}\label{6.28}
\hat{u}(t)=\frac{\tilde{\Lambda}(t)}{\Lambda(t)}(\hat{X}(t)+\phi^{-1}(t,\theta({t-}),y)\psi(t,\theta({t-}),y)),
\end{eqnarray}
where
\begin{eqnarray}\label{6.29}
\tilde{\Lambda}(t)={-\bar{m}(t,\theta({t-}))\sigma(t,\theta({t-}))+\int_{\Gamma}{g(t,X(t),\theta({t-}),\gamma)}\pi(d\gamma)}. \no \\
\Lambda(t)={\sigma}^{2}(t,\theta({t-}))+\phi(t,\theta({t-}),Y(t))\int_{\Gamma}{g^{'}(t,X(t),\theta({t-}),\gamma)g(t,X(t),\theta({t-}),\gamma)}\pi(d\gamma).
\end{eqnarray}
To find the optimal control it remains to find $\phi$ and $\psi$. To do  so set $X(t):= \hat{X}(t), u(t):=\hat{u}(t)$ and $p(t):=\hat{p}(t)$ in (\ref{6.22}) and then substitute for $\hat{p}(t)$ in (\ref{6.18}) and   $\hat{u}(t)$ from (\ref{6.28}) . As this result is linear in $\hat{X}(t)$ we compare the  coefficient on both side of the resulting equation to get following two equations namely,
\begin{eqnarray}\label{6.30}
0&=&2r\phi(t,i,Y(t))+\phi_{t}(t,i,Y(t))+\phi_{y}(t,i,Y(t))+\sum_{i \neq j,i=1}^{M}{p_{ij}\frac{f^{h}(y/i)}{1-F^{h}(y/i)}}{(\phi(t,j,0)-\phi(t,i,Y(t)))}\no\\
&+& \frac{\tilde{\Lambda}(t)}{\Lambda(t)}\sigma(t,i)\bar{m}(t,i)\phi(t,i,Y(t))-\frac{\tilde{\Lambda}(t)}{\Lambda(t)}\phi(t,i,Y(t))\int_{\Gamma}{g(t,X(t),i,\gamma)}\pi (d\gamma).
\end{eqnarray}
\begin{eqnarray}\label{6.31}
0&=&r\psi(t,i,Y(t))+\psi_{t}(t,i,Y(t))+\psi_{y}(t,i,Y(t))+\sum_{ i \neq j,i=1}^{M}{p_{ij}\frac{f^{h}(y/i)}{1-F^{h}(y/i)}(\psi(t,j,0)-\psi(t,i,Y(t)))}\no\\
&+&\frac{\tilde{\Lambda}(t)}{\Lambda(t)}\sigma(t,i)\bar{m}(t,i)\psi(t,i,Y(t))-\frac{\tilde{\Lambda}(t)}{\Lambda(t)}\psi(t,i,y)\int_{\Gamma}g(t,X(t),i,\gamma) \pi d(\gamma).\no\\
\end{eqnarray}
with terminal boundary conditions given by (\ref{6.19}). Consider the following process
\begin{eqnarray}\label{6.32}
\tilde{\phi}(t,i,y)=-2E\bigg[\exp\bigg\{\int_{t}^{T}\bigg(2r(s,\theta({s-}))+\frac{\tilde{\Lambda}(s)}{\Lambda(s)}\sigma(s,\theta({s-}))\bar{m}(s,\theta({s-}))\no\\
-\frac{\tilde{\Lambda}(s)}{\Lambda(s)}\int_{\Gamma}{g(s,X(s),\theta({s-}),\gamma)\pi(d\gamma)}\bigg)ds\bigg\}|{(\theta(s-)=i,Y(t)=y)}\bigg].\no\\
\end{eqnarray}
\begin{eqnarray}\label{6.33}
\tilde{\psi}(t,i,y)&=&2dE\bigg[\exp\bigg\{\int_{t}^{T}\bigg(r(\theta({s-}),s)+\frac{\tilde{\Lambda}(s)}{\Lambda(s)}\sigma(s,\theta({s-}))\bar{m}(s,\theta({s-}))\no\\
&-&\frac{\tilde{\Lambda}(s)}{\Lambda(s)}\int_{\Gamma}{g(s,X(s),\theta({s-}),\gamma)\pi(d\gamma)}\bigg)ds\bigg\}\bigg|{(\theta(s-)=i,Y(s)=y)}\bigg].\no\\
\end{eqnarray}
We aim to show that $\phi=\tilde{\phi}$  and $\psi=\tilde{\psi}$. We define the following martingales:
\begin{eqnarray}\label{6.34}
R(t)=
\resizebox{.9\hsize}{!}{$E\bigg[\exp\bigg\{\int_{0}^{T}\bigg(2r(s,\theta(s-))+\frac{\tilde{\Lambda}(s)}{\Lambda(s)}\sigma(s,\theta(s-))\bar{m}(s,\theta(s-))
-\frac{\tilde{\Lambda}(s)}{\Lambda(s)}\int_{\Gamma}{g(s,X(s),\theta(s-),\gamma)\pi(d\gamma)}\bigg)ds\bigg\}|\mathcal{F}_{t}^{\theta,y}\bigg]$},\no\\
\end{eqnarray}
\begin{eqnarray}\label{6.35}
S(t)=
\resizebox{.9\hsize}{!}{$E\bigg[\exp\bigg\{\int_{0}^{T}\bigg(r(s,\theta(s-))+\frac{\tilde{\Lambda}(s)}{\Lambda(s)}\sigma(s,\theta(s-))\bar{m}(s,\theta(s-))
-\frac{\tilde{\Lambda}(s)}{\Lambda(s)}\int_{\Gamma}{g(s,X(s),\theta(s-),\gamma)}\pi(d\gamma)\bigg)ds\bigg\}|\mathcal{F}_{t}^{\theta,y}\bigg]$},\no\\
\end{eqnarray}
where $\mathcal{F}_{t}^{\theta,y}$ is defined as usual. We follow steps similar to that as seen in Example 1 and conclude that $\phi=\tilde{\phi}$  and $\psi=\tilde{\psi}$ by using joint-Markov property of ($\theta(t),Y(t)$),  to obtain the following expression for the  control  $\hat{u}(t)$  given as
\begin{eqnarray*}
\hat{u}(t)=
\resizebox{.9\hsize}{!}{$\frac{\tilde{\Lambda}(t)}{\Lambda(t)}\bigg(\hat{X}(t)
-\frac{d E\bigg[\exp\bigg\{\int_{t}^{T}(r(s,\theta(s-))+\frac{\tilde{\Lambda}(s)}{\Lambda(s)}\sigma(s,\theta(s-))\bar{m}(s,\theta(s-))
-\frac{\tilde{\Lambda}(s)}{\Lambda(s)}\int_{\Gamma}{g(s,X(s),\theta(s-),\gamma)\pi(d\gamma)})ds\bigg\}|{(\theta(t-)=i,Y(t)=y)}\bigg]}{E\bigg[\exp\bigg\{\int_{t}^{T}(2r(s,\theta(s-))+\frac{\tilde{\Lambda}(s)}{\Lambda(s)}\sigma(s,\theta(s-))\bar{m}(s,\theta(s-))
-\frac{\tilde{\Lambda}(s)}{\Lambda(s)}\int_{\Gamma}{g(s,X(s),\theta(s-),\gamma)\pi(d\gamma)})ds\bigg\}|{(\theta(t)=i,Y(t)=y)}\bigg]}\bigg)$}.
\end{eqnarray*}
For the choice of the control parameter and the boundedness conditions on the market parameters $r,b$,$\sigma$ and $g$, the conditions of Theorem 4.1 are satisfied and hence $\hat{u}$ is the optimal control process.


\begin{thebibliography}{99}
\bibitem{CK}Cadenillas, A., and Karatzas, I. 1995. The stochastic maximum principle for linear, convex optimal
control with random coefficients. \textit{SIAM J. Control Optim.} 33(2): 590-–624.
\bibitem{CE}Cohen, S.N., Elliott, R.J. 2010. Comparisons for backward stochastic differential equations
on Markov chains and related no-arbitrage conditions. \textit{ Ann. Appl. Probab.} 20(1): 267-–311.
\bibitem{DG}Deshpande, A., and Ghosh M.K. 2008. Risk minimizing option pricing in a Markov modulated market. \textit{ Stochastic Anal.Appl.} 26(1): 313--324.
\bibitem{Do}Donnelly, C. 2011. Sufficient Stochastic Maximum Principle in a Regime-Switching Diffusion Model. \textit{ Appl Math Optim.} 64(2):155–-169.
\bibitem{Fr}Framstad N.C., Oksendal,B., and Sulem, A. 2004. Sufficient stochastic maximum principle for
the optimal control of jump diffusions and applications to finance. \textit{ J. Optim. Theory Appl.}121(1): 77–-98.
\bibitem{GG}Ghosh, M.K., and Goswami, A. 2009. Risk minimizing option pricing in a semi-Markov modulated market. \textit{ SIAM J. Control Optim.} 48(3): 1519--1541.
\bibitem{GS}Gikhman, I.I., and Skorokhod, A.V. 1983. \textit{ The Theory of Stochastic Processes II}, Springer, Berlin.
\bibitem{IW}Ikeda, N., and Watanabe, S. 1981. \textit{ Stochastic differential equations and diffusion processes}, North-Holland, Amsterdam.
\bibitem{KS}Karatzas, I., and  Shreve, S. 1998. \textit{Methods of Mathematical Finance}, Springer, New York.
\bibitem{Oks}Oksendal, B., and Sulem,A. 2007. \textit{Applied Stochastic Control of Jump Diffusions}, Second Edition, Universitext, Springer Verlag,  Berlin, Heidelberg.
    \bibitem{Zh}Zhang,X., Elliott,R.J., and  Siu, T.K. 2012.{ A Stochatic maximum principle for a Markov regime switching jump-diffusion model and its application to finance}. \textit{ SIAM J. Control Optim.} 50 (2): 964–-990.
\bibitem{Wh}Whittle,P. 1980. \textit{Risk sensitive Optimal control}, John Wiley and Sons.
\bibitem{XY} Zhou, X.Y., and  Yin G. 2003. Markowitz's mean-variance portfolio selection with regime switching:
A continuous-time model. \textit{ SIAM J. Control Optim.} 42(4): 1466-–1482.
\end{thebibliography}
\end{document}